\input amstex
\input xy
\input epsf
\xyoption{all}
\documentstyle{amsppt}
\document
\magnification=1200
\NoBlackBoxes
\nologo
\hoffset1.5cm
\voffset2cm
\def\D{\Cal{D}}
\def\G{\bold{G}}
\def\Q{\bold{Q}}
\def\P{\bold{P}}
\def\Z{\bold{Z}}
\def\a{\bold{a}}
\def\b{\bold{b}}
\def\e{\varepsilon}
\def\g{\gamma}
\pageheight {16cm}


\medskip

\centerline{\bf NON--COMMUTATIVE GENERALIZED DEDEKIND SYMBOLS} 
\bigskip

\centerline{\bf Yuri I. Manin}

\medskip

\centerline{Max--Planck--Institut f\"ur Mathematik, Bonn, Germany}

\bigskip

\hfill{\it To the memory of Andrey Todorov}

\bigskip

{\it ABSTRACT.}  We define and study generalized Dedekind symbols with values
in non--necessarily commutative groups, generalizing constructions of Sh. Fukuhara
in [Fu1], [Fu2].  Basic examples of such symbols are obtained by replacing
period integrals of modular forms (cf. [Ma1], [Ma2], [Kn1], [Kn2], [ChZ]) by iterated
period integrals introduced and studied in [Ma3], [Ma4].
\bigskip

\centerline{\bf 0. Introduction: classical and generalized Dedekind symbols}

\medskip

Let $\eta (z)=e^{\pi i z/12} \prod_{n=1}^{\infty} (1-e^{2\pi inz})$ be the Dedekind eta function.
From its modularity properties, it follows that for any $a\in PSL (2,\bold{Z})$,
the difference  $\roman{log}\,\eta(az)  - \roman{log}\,\eta(z)$ is an easy calculable function of $a$ and $z$.
After a normalization it turns out that the interesting part of this difference $d(p,q)$ (in the notation of  [Fu1])
depends only on the first column  $(p,q)$ of $a$, and satisfies functional equations
$$
d(p,q)=d(p,q+p), \quad d(p,-q)= -d(p,q),
\eqno(0.1)
$$
$$
d(p,q)-d(q,-p) = \frac{p^2+q^2-3pq+1}{12pq}.
\eqno(0.2)
$$
In the terminology of [Fu1], this means that $d(p,q)$ is the classical (odd)  Dedekind symbol
with reciprocity function\   $\dfrac{p^2+q^2-3pq+1}{12pq}$. Many beautiful applications
of Dedekind symbols were described in the Introduction to [KiMel].

\smallskip
In several papers Sh.~Fukuhara, following in particular T.~M.~Apostol ([Ap]),
 studied {\it generalized} Dedekind symbols. He defined them as functions $D(p,q)$
of couples of co--prime integers, satisfying equations (0.1) (with one possible sign change)
and equations (0.2) in which the right hand side could now be replaced by another function $F(p,q)$,
called {\it the reciprocity function} of $D$.

\smallskip

In particular, Fukuhara established functional equations for possible reciprocity functions and
gave a general method allowing one to reconstruct a generalized Dedekind symbol from
its reciprocity function.

\smallskip

Finally, replacing Dedekind eta functions by (period integrals of) other modular fiorms, Fukuhara calculated
many new reciprocity functions and classified them.

\smallskip

In this note we replace {\it commutative} value groups for $D$ and $F$ by {\it non--commutative} ones.

\smallskip

The section  1 is dedicated to the statement and proof of the non--commutative versions of  functional
equations (0.1), (0.2), and of Fukuhara's reconstruction procedure. 

\smallskip

In the section 2, we show how classes of such non--commutative Dedekind symbols and reciprocity
functions emerge very naturally if one replaces the
period integrals of modular forms by {\it iterated period integrals}, introduced and studied 
by the author in [Ma3], [Ma4].

\smallskip

Finally, section 3 clarifies the cohomological interpretation of Dedekind symbols.
Essentially, their reciprocity functions emerge as  components of certain non--commutative 1--cocycles
of the modular group  $PSL(2,\Z )$.

\medskip

{\it Acknowledgement.} I am grateful to Professor Shinji Fukuhara
whose work inspired this paper and who kindly sent me his remarks about  its
preliminary version.

\bigskip
\centerline{\bf 1.  Non--commutative Dedekind symbols and reciprocity functions}

\medskip

{\bf 1.1. Notation.} In the following, we choose a group $\G$, possibly non--commutative,
and write multiplicatively the composition law in it; $1_{\G}\in \G$ denotes the identity.  We put also
$$
W:= \{ (p,q)\in\bold{Z}\times \bold{Z} \,|\, \roman{gcd}\, (p,q)=1\}.
$$

\medskip
{\bf 1.2. Definition.} {\it A $\bold{G}$--valued reciprocity function  is a map $f: W\to \G$ satisfying the  
following conditions}
$$
f(p,-q)=f(-p,q).
\eqno(1.1)
$$
$$
f(p,q)f(-q,p)=1_{\G}
\eqno(1.2)
$$
$$
f(p,p+q)f(p+q,q)=f(p,q).
\eqno(1.3)
$$
\medskip

Notice that applying (1.3) to $p=1, q=0$, we get $f(1,1)=1_{\bold{G}}.$ From (1.2)
we then get $f(-1,1)=1_{\bold{G}}.$
Moreover,  $f(-p,-q)=f(p,q)$
so that $f(p,q)$ depends only on $q/p$.

\medskip

{\bf 1.3. Definition.} {\it Let $f$ be a $\bold{G}$--valued reciprocity function.

\smallskip
 A generalized  $\G$--valued Dedekind symbol $\Cal{D}$ with reciprocity
 function $f$ is a map
$$
\Cal{D}:\,W\to \G:\ (p,q)\mapsto \Cal{D}(p,q)
$$
satisfying the following conditions: 
$$
  \D(p,q)=\D(p, q+p),
\eqno(1.4)
$$
$$
 \D(p,-q) =\D(-p,q),
\eqno(1.5)
$$
so that $\D(-p,-q)= \D(p,q)$.
Finally, 
$$
\D(p,q)\D(q,-p)^{-1}=f(p,q).
\eqno(1.6)
$$
}
\smallskip

{\bf 1.3.1. Remarks.} a) Clearly, $f$ can be uniquely reconstructed from $\D$: this follows from (1.6).
Moreover, if $\D$ satisfies (1.4)--(1.6), then $f$ satisfies (1.1)--(1.3).

\smallskip

In fact, (1.1) follows from(1.6) and (1.5): 
$$
f(p,-q)= \D(p,-q)\D(-q,-p)^{-1} =  \D(-p,q)\D(q,p)^{-1} = f(-p,q).
$$
Similarly, (1.2) follows from (1.6) and (1.5):
$$
f(p,q)f(-q,p) = \D(p,q)\D(q,-p)^{-1} \D(-q,p)\D(p,q)^{-1}= 1_{\bold{G}}.
$$
Finally, to check (1.3), we get from (1.6):
$$
f(p,p+q) f(p+q,q)=   \D(p,p+q)\D(p+q,-p)^{-1} \D(p+q,q)\D(q, -p-q)^{-1}.
$$
Two middle terms cancel due to (1.4). The remaining two terms give
$$
\D(p,p+q) \D(q, -p-q)^{-1}  =  \D(p,q)\D(q,-p)^{-1}=f(p,q),
$$
again in view (1.4) and (1.6).

\smallskip

b) Let $f$ be a reciprocity function. Define a new function $g(p,q):= f(-p,q)$. A direct check shows 
that $g$ is a reciprocity function as well.  We may call $f$ {\it an even} reciprocity function,
if $g=f$ identically. If the group $\bold{G}$ is abelian, another  condition $g=f^{-1}$ defines
{\it odd} reciprocity function, and  each reciprocity function is the product of its even and odd components
(at least, if $\bold{G}$ is 2--divisible).
The basic reason for this is the fact that pointwise multiplication induces on the set 
of reciprocity functions a structure of abelian group.
This is the case in [Fu1], [Fu2] (where $\bold{G}$ is written additively). 

\smallskip

If $\bold{G}$ is non--abelian, however, then the oddity condition includes $f(p,q)=f(q,p),$
and, moreover, there is no obvious way to generate all reciprocity functions from even and odd ones.

\medskip

We will show below  that nevertheless, Fukuhara's method for reconstruction of a Dedekind symbol from
its reciprocity function works (with necessary modifications) not only for abelian groups of values, but also
for non--abelian ones. 

\medskip

{\bf 1.4. (Modified) continued fractions.}  Let $A_0,A_1, \dots ,A_n,\dots $
be commuting independent variables. Define by induction
polynomials $Q^{(n)}, P^{(n)}\in \bold{Z}[A_0,\dots , A_n]$:
$$
Q^{(0)}(A_0)= A_0,\ P^{(0)}(A_0)= 1,
$$
$$
Q^{(n+1)}(A_0,\dots ,A_{n+1})= A_0 Q^{(n)}(A_1,\dots ,A_{n+1}) - P^{(n)}(A_1,\dots ,A_{n+1}),
$$
$$
P^{(n+1)}(A_0,\dots ,A_{n+1}) = Q^{(n)}(A_1,\dots ,A_{n+1})
\eqno(1.7)
$$
One easily sees that 
$$
\frac{Q^{(n)}(A_0,\dots ,A_{n})}{P^{(n)}(A_0,\dots ,A_{n})}=
\langle A_0,A_1,\dots A_n \rangle := A_0-\frac{1}{A_1-\frac{1}{A_2- \dots}}
\eqno(1.8)
$$
where the lowest layer  in the continued fraction representation is $A_{n-1}-\dfrac{1}{A_n}$.
For any integer values of $A_i$ and each $n$, values of $P^{(n)}, Q^{(n)}$ are coprime.


\medskip

{\bf 1.5. Basic moves.} Consider finite sequences of integers $\a=(a_0,\dots ,a_n)$, with arbitrary $n\ge 0$.
Define three types of moves on this set (cf. [Fu1], p. 88):

\smallskip

(i) This move depends on a choice of $\e \in \{1,-1\}$ and $i\in \{0,\dots ,n-1\}$
and makes the sequence longer by one:
$$
\a \mapsto \b:= (a_0,\dots , a_i+\e, \e, a_{i+1}+\e, a_{i+2} \dots ).
\eqno(1.9)
$$

(ii) A move of the second type depends on a choice of  $i\in \{0,\dots ,n\}$ and a representation
of $a_i$ as $a_i=b+c$. It 
makes the sequence longer by two:
$$
\a \mapsto \b^{\prime}:= (a_0,\dots , a_{i-1}, b,0,c, a_{i+1} , \dots , a_n).
\eqno(1.10)
$$

(iii) Finally, a move of the third type  depends on a choice of $\e \in \{1,-1\}$ 
and makes the sequence longer by one:
$$
\a \mapsto \b^{\prime\prime}:= (a_0,\dots , a_{n-1}, a_{n}+\e, \e).
\eqno(1.11)
$$

{\bf 1.6.1. Lemma.} {\it (i) 
If two sequences $\a$  differ by  one of the moves (1.8)--(1.10),
the respective continued fractions define one and the same rational  number.

\smallskip

(ii) Conversely, if a rational number is written in two different ways as a
continued fraction (1.11), then the respective sequences differ by a finite sequence
of moves (1.8)--(1.10) and their inverses.}

\smallskip

The first part is straightforward. For the second one, see [Mel] and [KiMel].

\medskip

{\bf 1.7. Explicit presentation of the  symbol.} We will now show that Fukuhara's method  in [Fu1]
for reconstruction of a generalized Dedekind symbol from its reciprocity function can be applied
also in non--commutative case.
\smallskip

For a given $(p,q)\in W$, choose a presentation as above
$\dfrac{q}{p}=\langle a_0, \dots ,a_n\rangle$ that is, more precisely,
$$
q=q_0:=Q^{(n)}(a_0,\dots ,a_n),\quad   p=p_0:=P^{(n)}(a_0,\dots ,a_n)
\eqno(1.12)
$$

Furthermore, for $i=1, \dots , n$, define $(q_i,p_i)\in \bold{Z}^2$  by
$$
q_i:=Q^{(n-i)}(a_i,\dots ,a_n),\   p_i:=P^{(n-i)}(a_i,\dots ,a_n)
\eqno(1.13)
$$
(cf. (1.7), (1.8)) so that  $\dfrac{q_i}{p_i}=\langle a_i, \dots ,a_n\rangle$ .
\smallskip
Finally, put 
$$
\D(p,q):=f(p_1,q_1)^{-1} f(p_2,q_2)^{-1} \dots f(p_n,q_n)^{-1} \in \G.
\eqno(1.14)
$$
\medskip

{\bf 1.8. Theorem.} {\it (i) $\D(p,q)$ depends only on $(p,q)$ and not on the choice of the
presentation  (1.14).

(ii) The function $\D$ thus defined is a generalized Dedekind symbol
with reciprocity function  $f$.}

\smallskip

{\bf Proof.} (i) In view of Lemma 1.6.1, it suffices to check that the r.h.s. of  (1.14) does not change whenever
one applies to $\a:= (a_0, \dots ,a_n)$ one of the moves (1.9)--(1.11)

\smallskip

{\it MOVE (1.9).} Let $\b=(b_0,b_1, \dots, b_{n+1})$ and $t_j/s_j:= \langle b_j,\dots ,b_{n+1} \rangle.$
Then for $j\le i$ we have $t_j/s_j=q_j/p_j$ and for $j\ge i+3$ we have  $t_j/s_j=q_{j-1}/p_{j-1}$. 
Hence two products in the rhs of (1.14) corresponding to $\a$ and $\b$ respectively
have the same beginnings and endings, and we have to check only the coincidence
of the respective middle segments. Passing to their inverses, we will have to prove that
$$
  (?)\quad\quad  f(p_{i+1}, q_{i+1})= f(s_{i+2},t_{i+2})f(s_{i+1},t_{i+1}).
 \eqno(1.15)
 $$ 
 
From (1.7) and (1.13), we get the identities:
 $$
 p_{i+1}=q_{i+2},\  q_{i+1}=a_{i+1}q_{i+2}-p_{i+2};
 $$
 $$
 s_{i+2}=q_{i+2},\  t_{i+2}= (a_{i+1}+\e)q_{i+2}-p_{i+2};
 $$
 $$
 s_{i+1}= (a_{i+1}+\e)q_{i+2}-p_{i+2},\  t_{i+1}= \e(a_{i+1}q_{i+2}-p_{i+2}).
 $$
Hence, comparing with (1.14), we see that (1.15) is equivalent to the identity
$$
f(q_{i+2}, a_{i+1}q_{i+2}-p_{i+2}) 
$$
$$
=f(q_{i+2}, (a_{i+1}+\e)q_{i+2}-p_{i+2}) f( (a_{i+1}+\e)q_{i+2}-p_{i+2},\e(a_{i+1}q_{i+2}-p_{i+2}) ).
\eqno(1.16)
$$
To make our check of (1.16) more transparent, denote the three expressions $f(...)$ in (1.16)
respectively $A,B,C$ so that (1.16) means $A=BC$. Moreover, denote the subsequent pairs
of arguments in $A,B,C$ by  $(1), (2); (3), (4);  (5),(6)$ respectively so that e.g.  $(3)=q_{i+2}$.

\smallskip

Now the check works slightly differently in the following two cases:

\smallskip

(a) $\e=1.$ Then $(4)=(5)=(3)+(6)$, and we get (1.16) after substituting  $p\mapsto (3)=(1)$
and $q\mapsto (6)=(2)$ into (1.3).

\smallskip

(b) $\e=-1.$ Then the same strategy works for the identity $AC^{-1}=B$ which is equivalent to 
$A=BC$. In fact, from (1.2) we get
$$
C^{-1}= f(a_{i+1}q_{i+2}-p_{i+2}, (a_{i+1}-1)q_{i+2}-p_{i+2}).
$$
so that $AC^{-1}=B$ becomes
$$
f(q_{i+2}, a_{i+1}q_{i+2}-p_{i+2}) f(a_{i+1}q_{i+2}-p_{i+2}, (a_{i+1}-1)q_{i+2}-p_{i+2}) 
$$
$$
=f(q_{i+2}, (a_{i+1}-1)q_{i+2}-p_{i+2})
$$
which is again a special case of (1.3).

\medskip

{\it MOVE (1.10).} Let $\b^{\prime}=(b_0,b_1, \dots, b_{n+2})$ and $t_i^{\prime}/s_i^{\prime}:= \langle b_i,\dots ,b_{n+2} \rangle.$
Then for $j\le i$ we have $t_j^{\prime}/s_j^{\prime}=q_j/p_j$ and for $j\ge i+3$ we have  $t_j^{\prime}/s_j^{\prime}=q_{j-2}/p_{j-2}$. 
Hence the two products in the r.h.s. of (1.14) corresponding to $\a$ and $\b$ respectively
have the same beginnings and endings, but passage to $\b^{\prime}$ introduces extra terms in--between:
$f(s^{\prime}_{i+1},t^{\prime}_{i+1})^{-1}$ and $ f(s^{\prime}_{i+2},t^{\prime}_{i+2})^{-1}$. In our setting, their product (in any order)
is  $1_{\bold{G}}$
thanks to  (1.2),  because
a straightforward check shows that
$$
\frac{t^{\prime}_{i+1}}{s^{\prime}_{i+1}}=\frac{-q_{i+1}}{cq_{i+1}-p_{i+1}},\quad
\frac{t^{\prime}_{i+2}}{s^{\prime}_{i+2}}=\frac{cq_{i+1}-p_{i+1}}{q_{i+1}}.
$$

\medskip

{\it MOVE (1.11).} Let $\b^{\prime\prime}=(b_0,b_1, \dots, b_{n+1})$ and 
$t_i^{\prime\prime}/s_i^{\prime\prime}:= \langle b_i,\dots ,b_{n+1} \rangle.$
Then for $j\le n$ we have $t_j^{\prime\prime}/s_j^{\prime\prime}=q_j/p_j$, and moreover  
$t_{n+1}^{\prime\prime}/s_{n+1}^{\prime\prime}=\e/1$.  Hence this extra term contributes
$f(\e,1)=1_{\bold{G}}$.

\medskip

(ii) We will now check that  products (1.14) satisfy  identities (1.4)--(1.6).

\smallskip

In fact, if $q/p=\langle a_0,a_1, \dots ,a_n\rangle$, then $(p+q)/p =  \langle a_0+1,a_1, \dots ,a_n\rangle$,
so that  $(1.14)$ gives for $\D(p,q)$  and $\D(p,p+q)$ the same expressions. This shows (1.4).

\smallskip

Furthermore, in the same notation we have $-q/p=\langle -a_0,-a_1, \dots ,-a_n\rangle$,
$-q_i/p_i=\langle -a_i, \dots ,-a_n\rangle$,  so that 
(1.14) gives the same expressions for $\D(p,-q)$ and  $\D(-p,q)$.  This shows $(1.5)$.

\medskip

 It remains to establish $(1.6)$. If  $q/p=\langle a_0,a_1,\dots ,a_n\rangle$, then
 $-p/q=\langle 0, a_0, a_1, \dots ,a_n\rangle .$
 From (1.14) one sees that 
 $$
\D(p, q)\D(q,-p)^{-1}=   f(p_0,q_0) =f(p,q)
$$
From (1.2) it follows that this is equivalent to (1.6).

\newpage
\centerline{\bf 2.  Reciprocity functions from iterated integrals}

\medskip

{\bf 2.1. Notations.} Let $w\ge 0$ be an even integer. Denote by $\{\varphi_j(z)\}, j=1,\dots r$,
a basis of the space of  cusp forms of weight $w+2$ for the full modular group
$PSL(2,\bold{Z}).$ We have $\varphi_j(z+1)=\varphi_j(z)$
and $\varphi_j(-z^{-1})=  \varphi_j(z)z^{w+2}$ for all $j$.
\smallskip

Let $(A_j),  j=1,\dots ,r$, be independent associative but non--commuting formal variables.
For $(p,q)\in W$, put
$$
\Omega (p,q):=\sum_{j=1}^r A_j \varphi_j(z) (pz-q)^w\,dz,
\eqno(2.1)
$$

Finally, in the notations of [Ma3], sec.~1, consider iterated integrals along the geodesics from rational points to 
$i\infty$ in the upper half plane:
$$
f(p,q) := J_0^{i\infty}(\Omega (p,q)), \quad \D (p,q):=  J_{p/q}^{i\infty}(\Omega (p,q))
\eqno(2.2)
$$
taking values in the ring of free associative formal series $\bold{C}\langle\langle A_1,\dots ,A_r\rangle\rangle$,
or, more precisely, in the multiplicative subgroup
$\bold{G}:= 1+(A_1,\dots ,A_r)$ of this ring.

\smallskip

We remind briefly that all  iterated integrals that we will be considering below are taken along oriented geodesics
in the upper complex half--plane with coordinate $z$, connecting two cusp points $a,b \in \bold{P}^1(\bold{Q})$:
rational points of the real $z$--line and infinity, and we use two types of relations among such integrals:
$$
J^{a_1}_{a_2} (\Omega)J^{a_2}_{a_3} (\Omega) \dots  J_{a_n}^{a_{n-1}} (\Omega)J^{a_n}_{a_1} (\Omega)=1
$$
(cf. [Ma3], (1.9)), and variable change, or functoriality (cf. [Ma3], (1.10)).
\smallskip

Notice that our group $\bold{G}$  becomes non--commutative only for $r>1$ that is, for $w> 26$.
One can get examples of lesser weight, generalizing the constructions below to congruence subgroups and/or
Eisenstein series. For constructions with Eisenstein series (regularization of iterated integrals) cf.
section $6$ of [Ma3]. 

\medskip

{\bf 2.2. Theorem.} {\it The map $ \D$ : $W\to\bold{G}$ is  a generalized Dedekind symbol with reciprocity function $f$.}

\smallskip

{\it Remark.} In the preliminary version of this paper,  only formula for the reciprocity function $f$
was given explicitly, together with a proof that it satisfies all necessary identities and therefore 
that $\D(p,q)$ can be in principle reconstructed from it. Professor Fukuhara found  the simple explicit formula
for $\D$ given here and kindly allowed me to include it here  (e-message of 13/11/2012).

\smallskip

In fact, the linear in $(A_i)$ term of $\D(p,q)$ is precisely the linear combination of
commutative generalized Dedekind symbols from [Fu1], Definition 7.1.

\smallskip

{\bf Proof.} We have to prove relations (1.1)--(1.3) for $f$.

\smallskip

 {\it Relations (1.1).} Since $w$ is even, we have
$$
(pz+q)^w=(-pz-q)^w.
$$
This shows that   $f(p,-q)=f(-p,q).$

\smallskip

{\it Relations (1.2).} In the iterated integral   $f(p,q) := J_0^{i\infty}(\Omega (p,q))$
we can make the variable change $z=-u^{-1}$. Then from (2.1) and $\varphi_j(-u^{-1})=\varphi_j(u)u^{w+2}$ we see that
$$
\Omega (p,q)=\sum_{j=1}^r A_j \varphi_j (-u^{-1}) (-pu^{-1}-q)^w\,d(-u^{-1})=
\sum_{j=1}^r A_j \varphi_j (u) (-qu-p)^w\,du = \Omega (-q,p).
$$
Moreover, the integration limits in the variable $u$ will be $(i\infty ,0)$. Hence finally
$$
f(p,q) = J_0^{i\infty}(\Omega (p,q)) = J^0_{i\infty}(\Omega(-q,p)) =f(-q,p)^{-1}.
$$

\smallskip

{\it Relations (1.3).}  We must check the identity
$$
J_0^{i\infty}(\Omega (p,p+q)) J_0^{i\infty}(\Omega (p+q,q))=
J_0^{i\infty}(\Omega (p,q))
$$
or equivalently,
$$
J_0^{i\infty}(\Omega (p+q,q))=
J^0_{i\infty}(\Omega (p,p+q)) J_0^{i\infty}(\Omega (p,q))
$$
We rewrite the last three integrals in turn.

\smallskip

The l.h.s.  integral  becomes
$J^0_{i\infty}(\Omega (-q, p+q))$ after the variable change $z=-u^{-1}$.
\smallskip

The first r.h.s. integral becomes 
$J^{-1}_{i\infty}(\Omega (p,q))$ after the variable change $z=u+1$.

\smallskip
Hence the product of the two r.h.s. integrals equals $J^{-1}_{0}(\Omega (p,q))$. Let us make in this last integral
the substitution $z=-(u+1)^{-1}$. The integration limits $(0,-1)$ w.r.t $z$ will become $( i\infty ,0 )$ w.r.t $u$. We have

\medskip
$$
\varphi_j(-(u+1)^{-1}) \,(-p(u+1)^{-1}-q)^w\, (u+1)^{-2}du
$$
$$
=\varphi_j(u)(p+qu+q)^{w}\,du.
$$
 Hence finally the r.h.s integral becomes  $J^0_{i\infty}(\Omega (q, -p-q))$ that coincides with
$J^0_{i\infty}(\Omega (-q, p+q))$ in view of (1.1).

\smallskip

It remains to check the relations (1.4)--(1.6). The first one is obtained by the variable change $z=u+1$, the second is obvious. Finally,
(1.6) is obtained  by the substitution $z=-u^{-1}$ showing that
that  $\D(q,-p)^{-1}=J_0^{q/p}(\Omega(p,q))$. This completes the proof.

\medskip

{\it Remark.} The identities (1.1)--(1.3) for our $f$ form an iterated version of the classical Shimura--Eichler 
relations for periods: cf. [Ma3], Proposition 2.1. We have essentially reproduced here the structure
of the classical proof.  Since these computations have a well known cohomological interpretation,
we may expect that  our generaliized Dedekind symbols and their reciprocity functions also can be interpreted in this way.
We develop this remark in the next section.
\bigskip

\centerline{\bf 3. Reciprocity functions as non--commutative cocycles}

\medskip

{\bf 3.1. Non--commutative 1--cohomology.} Below we consider abstract groups, generally non--commutative,
and group laws are written multiplicatively. 
 Let $\Gamma$ be a group, and $\G$ a group endowed with a left action of $\Gamma$ by automorphisms: $(\gamma,g)\mapsto 
 \gamma g$.

\smallskip

We define 1--cocycles by
$$
Z^1(\Gamma ,\G) := \{\, u:\,\Gamma \to \G\,|\,u(\gamma_1\gamma_2)= u(\g_1)\,\g_1u(\g_2)\,\} .
\eqno(3.1)
$$
It follows that $u(1_{\Gamma})=1_{\G_0}$.

\smallskip

Two cocycles are equivalent, $u^{\prime} \sim u$, iff
there exists a $g\in \G$ such that for all $\g \in \Gamma$ we have $u^{\prime}(\g)=g^{-1}\,u(\g)\cdot \g g.$ This is an equivalence relation,
and by definition, 
$$
H^1(\Gamma ,\G) := Z^1(\Gamma ,\G)/(\sim ).
$$
This is a set with a marked point: the class of the
trivial cocycles  $u_{g}(\g) =g^{-1} \g g$.

\medskip

{\bf 3.2. Cohomology of $PSL(2,\bold{Z})$.} Consider now the group
$\Gamma =PSL(2,\bold{Z})$, and let $\G$ be a possibly noncommutative $\Gamma$--module.
It is known that $PSL(2,\bold{Z})$ is the free product of
its two subgroups $\bold{Z}_2$ and $\bold{Z}_3$ generated respectively by
$$
\sigma = \left(\matrix 0& -1\\ 1& 0\endmatrix \right),\quad
\tau = \left(\matrix 0& -1\\ 1& -1\endmatrix \right)\, .
$$

$PSL(2,\bold{Z})$ acts transitively on $\bold{P}^1(\bold{Q}),$
the set of cusps of upper complex half--plane.
The stabilizer of $\infty$ is a cyclic subgroup $G_{\infty}$ generated by $\sigma\tau$. Hence the stabilizer $G_a$ of any cusp
$a\in \bold{P}^1(\bold{Q})$ is generated by $\g^{-1}\sigma\tau \g$
where $\g a =\infty$.

\smallskip

Following [Ma4], we will give a short description of the set $H^1(PSL(2,\bold{Z}),\G)$
and its cuspidal subset $H^1(PSL(2,\bold{Z}),\G)_{cusp}$
consisting by definition of those cocycle classes that become trivial
after restriction to any $G_a$.

\medskip

{\bf 3.3. Proposition.} {\it (i) Restriction to  $(\sigma ,\tau )$
of any cocycle in  $Z^1(PSL(2,\bold{Z}), \G)$
belongs to the set
$$
\{\, (X,Y) \in \G\times \G\,|\, X\cdot\sigma X=1,\,
Y\cdot\tau Y\cdot\tau^2Y=1\,\}.
\eqno(3.2)
$$

\smallskip

(ii) Conversely, any element of the set (3.2) comes from
a unique 1--cocycle so that we can and will identify
these two sets. The cohomology 
relation between cocycles translates as
$$
(X,Y) \sim (g^{-1}X\sigma g, g^{-1}Y\tau g), \  g\in \G.
\eqno(3.3)
$$

\smallskip

(iii) The cuspidal part of the cohomology consists of
classes of pairs of the form
$$
\{\,(X,Y)\,|\,\exists g\in \G,\, X\cdot\sigma Y=g^{-1}\cdot \sigma\tau g\,\}.
\eqno(3.4)
$$}
We may call (3.2) abstract (noncommutative) {\it Shimura--Eichler relations.} 

\medskip

{\bf 3.4. Definition.} {\it An element $(X,Y)$ of (3.2) is called (the representative of)  a Dedekind cocycle, iff
it satisfies the relation
$$
Y=\tau X. \eqno(3.5)
$$
}

\medskip

{\bf 3.5. Reciprocity functions as cocycles.} Let now $\G_0$ be a group. Denote by $\G$
the group of functions $f:\,   \P^1(\Q)\to \G_0$ with pointwise multiplication. Define the left
action of $\Gamma$ upon $\bold{G}$ by
$$
(\gamma f)(x) = f(\g^{-1}x);\quad f\in \G,\ x\in \P^1(\Q),\ \gamma \in \Gamma.
\eqno(3.6)
$$

Let $f:\, W\to \G_0$ be a $\G_0$--valued reciprocity function, as in Def. 1.2.
Define elements $X_f, Y_f\in \G$ as the following functions $\P^1(\Q)\to \G_0$:
$$
X_f(qp^{-1}):= f(p,q),
$$
$$
Y_f(qp^{-1}):= (\tau X_f)(qp^{-1})= X_f(\tau^{-1}(qp^{-1}))= f(q, q-p).
\eqno(3.7)
$$
\medskip

{\bf 3.6. Theorem.} {\it    The map $f\mapsto (X_f,Y_f)$ establishes a bijection between the
set of $\G_0$--valued reciprocity functions and the set of (representatives of)  Dedekind cocycles
from $Z^1(\Gamma ,\G ).$  }

\medskip

{\bf Proof.} First, we check that if $f$ is a reciprocity function, then the pair $(X_f,Y_f)$ satisfies (3.2). In fact,
$$
(X_f\cdot  \sigma X_f)(qp^{-1}) = X_f(qp^{-1}) X_f(-pq^{-1})= f(p,q)f(-q,p)=1_{\G_0}
$$
in view of (1.2). Furthermore, from (3.7) we find
$$
(Y_f\cdot  \tau Y_f  \cdot  \tau^2 Y_f)(qp^{-1})= f(q,q-p)f(q-p,-p)f(p,q) =  f(q,-p)f(p,q)=   1_{\G_0}
$$
in view of (1.3) and (1.2).

\smallskip

Conversely, let $(X, Y=\tau X)$ be a Dedekind cocycle. Define the function $f:\, W\to \G_0$
 by $f(p,q):= X(qp^{-1})$ so that (1.1) is straightforward. Substituting this into (3.2), we get 
from previous computations (1.2) and (1.3).

\bigskip
\centerline{\bf References}

\medskip

[Ap] T.~M.~Apostol. {\it Generalized Dedekind sums and transformation formulae of certain Lambert series.} Duke 
Math.~J. 17 (1950), 147--157.

\smallskip

[ChZ] Y.~J.~Choie, D.~Zagier. {\it Rational period functions
for $PSL(2,\bold{Z} ).$} In: A Tribute to Emil Grosswald:
Number Theory and relatied Analysis, Cont. Math.,
143 (1993), AMS, Providence, 89--108.

\smallskip

[Fu1] Sh.~Fukuhara. {\it Modular forms, generalized Dedekind symbols
and period polynomials.} Math. Ann., 310 (1998), 83--101.

\smallskip

[Fu2] Sh.~Fukuhara. {\it Dedekind symbols with polynomial reciprocity laws.}
Math. Ann., 329 (2004), No. 2,  315--334.

\smallskip

[KiMel] R.~C.~Kirby,  P.~Melvin. {\it Dedekind sums, $\mu$--invariants and the signature cocycle.} 
Math. Ann. 299 (1994), 231--267.

\smallskip

[Kn1] M.~Knopp. {\it Rational period functions of the modular group.}
Duke Math. Journal, 45 (1978), 47--62.

\smallskip

[Kn2] M.~Knopp. {\it Rational period functions of the modular group II.}
Glasgow Math. Journal, 22 (1981), 185--197.

\smallskip

[Ma1] Yu.~Manin. {\it Parabolic points and zeta-functions of modular curves.}
Russian: Izv. AN SSSR, ser. mat. 36:1 (1972), 19--66. English:
Math. USSR Izvestiya, publ. by AMS, vol. 6, No. 1 (1972), 19--64,
and Selected Papers, World Scientific, 1996, 202--247.

\smallskip

[Ma2] Yu.~Manin. {\it Periods of parabolic forms and $p$--adic Hecke series.}
Russian: Mat. Sbornik, 92:3 (1973), 378--401. English:
Math. USSR Sbornik, 21:3 (1973), 371--393,
and Selected Papers, World Scientific, 1996, 268--290.

\smallskip

[Ma3] Yu.~Manin. {\it Iterated integrals of modular forms and 
noncommutative modular symbols.}
In: Algebraic Geometry and Number Theory.
In honor of V. Drinfeld's 50th birthday.
Ed. V.~Ginzburg. Progress in Math., vol. 253. 
Birkh\"auser, Boston, pp. 565--597. Preprint math.NT/0502576.

\smallskip

[Ma4] Yu.~Manin. {\it Iterated Shimura integrals.} Moscow Math. Journal, 
vol. 5, Nr. 4 (2005), 869--881.
 Preprint math.AG/0507438.

\smallskip

[Mel] P.~Melvin. {\it Tori in the diffeomorphism groups of of simply--connected  4--manifolds.}
Math. Proc. Cambridge Phil. Soc. 91 (1982),  305--314.

\smallskip

\enddocument